\documentclass[titlepage]{amsart}
\newtheorem{prop}{Proposition}[section]
\newtheorem{lema}[prop]{Lemma}
\newtheorem{teo}[prop]{Theorem}
\newtheorem{corolario}[prop]{Corollary}

\newtheorem{remark}[prop]{\sc Remark}

\title{On the Dunford-Pettis property of the  tensor
product of $C(K) $ spaces }

\author{Fernando Bombal}
\address{Departamento de An lisis Matem tico \\
Facultad de Matem ticas \\ Universidad Complutense de Madrid \\
Madrid 28040 }
\email{bombal@eucmax.sim.ucm.es, ignacio\_villanueva@mat.ucm.es}

\author{Ignacio Villanueva }
 \thanks{Both authors were partially supported by DGICYT grant
PB97-0240.}

\begin{document}

\begin{abstract}
In this paper we characterize those compact Hausdorff spaces $K$
such that $C(K)\hat{\otimes}C(K)$ (and $C(K)\hat{\otimes}_s C(K)$)
have the Dunford-Pettis Property, answering thus in the negative a
question posed in \cite{CaGo} which asked if
$\ell_{\infty}\hat{\otimes}\ell_{\infty}$ and
$C[0,1]\hat{\otimes}C[0, 1]$ have this property.
\end{abstract}

\subjclass{46B28, 47B07}

\keywords{Dunford-Pettis property, spaces of continuous functions,
projective tensor product}
\maketitle

\section{Introduction}
In \cite{CaGo} it is asked if
$\ell_{\infty}\hat{\otimes}\ell_{\infty}$ and
$C[0,1]\hat{\otimes}C[0,1]$ have the Dunford-Pettis Property (DPP).
The purpose of this paper is to show that this is not the case, and
to characterize those $C(K) $ spaces such that $\hat{\otimes}_n
C(K)$ or $\hat{\otimes}_{n,s} C(K)$ have the Dunford-Pettis
Property, when $n>1$.

Let us recall that a Banach space $E$ has the DPP if, for any
Banach space $F$, every weakly compact operator from $E$ into $F$
is completely continuous. In this paper we will need the following
facts about the DPP; they are all well known, and can be found, for
instance, in \cite{Di}: (a) $C(K)$ spaces and Schur spaces have the
DPP, (b) if a dual Banach space $E^*$ has the DPP, then so does $E$
(but the reciprocal is not true), (c) the DPP is stable under
complemented subspaces.

We will now explain our notation. Throughout the paper $K$ will
denote a compact Hausdorff space and $C(K)$ the space of continuous
scalar functions defined on it. We will write $\mathcal L(E;X)$ to
indicate the linear operators from $E$ into $X$, and $\mathcal
L^n(E_1,\ldots, E_n
;X)$ will denote the space of continuous multilinear operators from
$E_1 \times \cdots \times E_n$ into $X$. When $E_1= \cdots = E_n
=E$ we will write this as  $\mathcal L^{n}(E;X)$ and in this case
$\mathcal L^{n}_s (E;X)$ will refer to those of the previous
operators which are symmetric. In all the cases, if $X=\Bbb K$, the
scalar field, we will not write it. We will write the projective
tensor product of $E_1,\ldots, E_n$ as $E_1\hat{\otimes}\cdots
\hat{\otimes} E_n$. We will write $\hat{\otimes}_n E$ to denote the
$n$-fold projective tensor product of $E$, and $\hat{\otimes}_{n,s}
E$ will denote the symmetric $n$-fold projective tensor product of
$E$. If $n=2$ we will sometimes also write this as
$E\hat{\otimes}_s E$.

We consider to be well known that $\mathcal L^n(E_1,\ldots, E_n
;X)$ is naturally isometric to $\mathcal L (E_1\hat{\otimes}\cdots
 \hat{\otimes} E_n;X)$ and that $\mathcal L^{n}_s(E;X)$ is naturally
isometric to $\mathcal L (\hat{\otimes}_{n,s} E ;X)$.

We will use the convention $\stackrel{[i]}{\ldots}$ to indicate
that the $i^{th}$ coordinate is not involved.

\section{The results}

As is well known, every linear operator from any $C(K)$ space into
the dual of another $C(K)$ space is weakly compact, and therefore
completely continuous. Using that we can prove the following lemma.

\begin{lema}\label{lema} Let $K_1$, $K_2 $ be two compact Hausdorff spaces.
Let $(f_n)\subset C(K_1)$ be a weakly null sequence and let
$(g_n)\subset C(K_2)$ be a bounded sequence. Then the sequence
$(f_n \otimes g_n)\subset C(K_1)\hat{\otimes}C(K_2)$ is weakly
null.
\end{lema}
\begin{proof}
Let $\phi \in (C(K_1)\hat{\otimes}C(K_2))^* $, and let us consider
the linear operator $S\in \mathcal L(C(K_1);C(K_2)^*) $ associated
to it defined by
     $$S(f)(g)=\phi(f\otimes g) \quad.$$
Let us suppose without loss of generality that $\sup_n \|g_n\| \leq
1$. Since $S$ is completely continuous, we get that
    $$\lim_{n\rightarrow\infty} |\phi(f_n\otimes g_n)| \leq
    \lim_{n\rightarrow\infty} \|S(f_n)\| =0 \quad.$$
\end{proof}

In \cite{BoVi} it can be seen that, if $K$ is scattered, then
$(\hat{\otimes}_n C(K))^*$ and $(\hat{\otimes}_{n,s} C(K))^*$ are
Schur spaces for every $n\in \Bbb N$. With analogous proof it can
be proved that, in case $K_1,\ldots K_n$ are scattered,
$(C(K_1)\hat{\otimes}\cdots\hat{\otimes}C(K_n))^*$ is a Schur
space.

We state now our main result.

\begin{teo}\label{teoprincipal}
Let $K_1 , K_2$ be two infinite compact Hausdorff spaces. Then
\linebreak $C(K_1)\hat{\otimes} C(K_2)$ has the DPP if and only if
both $K_1$ and $K_2$ are scattered.
\end{teo}

\begin{proof} If both $K_1$ and $K_2$ are scattered, then $(C(K_1)
\hat{\otimes}C(K_2))^*$ is a Schur space and therefore $C(K_1)
\hat{\otimes}C(K_2)$ has the DPP. Now, let us suppose that one of
them, say $K_2$, is not scattered. Since $K_1$ is infinite,
$C(K_1)$ is not Schur, and therefore there exist two sequences
$(f_n)\subset B_{C(K_1)}$ and $(\xi_n)\subset B_{(C(K_1))^*}$ such
that $(f_n)$ is weakly null and $\xi_n(f_n)=1$ for every $n\in \Bbb
N$. Also, since $K_2$ is not scattered, $C(K_2)$ contains an
isomorphic copy of $\ell_1$, and therefore there exists a
continuous surjective operator $q: C(K_2)\mapsto \ell_2$
(\cite[Corollary 4.16]{DiJT}). Then let
us consider the trilinear form
$$ T: C(K_1) \times C(K_2) \times C(K_2) \mapsto \Bbb K $$ defined
by
$$ T(f,g,h)=\sum_{n=1}^{\infty} \xi_n(f) q(g)_n q(h)_n $$
and let us consider the linear operator
    $$\hat{T}^1:C(K_1) \hat{\otimes} C(K_2) \mapsto
(C(K_2))^*$$ canonically associated to it given by
$$ \hat{T}^1(f\otimes g)(h)=T(f,g,h) \quad.$$

It is clear that $\hat{T}^1 = q^* \circ \psi \circ \phi $ where
    $$\phi:C(K_1) \hat{\otimes} C(K_2)\mapsto \ell_2 $$
is given by
    $$\phi(f\otimes g) = (\xi_n(f) q(g)_n)_n\quad.$$
and $\psi\in\mathcal L(\ell_2;\ell_2^*)  $ is the canonical linear
isometry identifying both spaces. Since $\psi$, $\phi$ and $q^*$
are all of them weakly compact, so is $\hat{T}^1$.

So we now just have to see that $\hat{T}^1$ is not completely
continuous. Let us consider a sequence of bounded functions
$(g_n)_n\subset C(K_2)$ such that $q(g_n)=e_n$, where $(e_n)$ is
the canonical basis of $\ell_2$. Then, according to Lemma
\ref{lema}, the sequence $(f_n \otimes g_n)_n\subset
C(K_1)\hat{\otimes}C(K_2)$ weakly converges to zero, but, for each
$n\in \Bbb N$,
    $$ \| \hat{T}^1(f_n \otimes g_n)\| \sup_n \|g_n\| \geq |\hat{T}^1
    (f_n \otimes g_n)(g_n)| = |T(f_n, g_n, g_n)| = 1$$
a contradiction.
\end{proof}

 \begin{remark}
{\rm In \cite{ArGa}, the following definition is stated: a
multilinear form $T\in \mathcal L^n(E_1, \ldots, E_n)$ is said to
be
\begin{em}regular\end{em} if every one of the associated linear operators
    $$T_1^i: E_i \mapsto \mathcal L^{n-1}(E_1,
    \stackrel{[i]}{\ldots}, E_n)$$
is weakly compact, and this is shown to be equivalent to every one
of the associated $(n-1)$-linear operators
    $$T_{n-1}^i: E_1 \times \stackrel{[i]}{\cdots} \times E_n
     \mapsto \mathcal L(E_i)$$
being weakly compact. From the proof it is clear that, given an
$i\in (1,\ldots, k)$, $T_1^i$ is weakly compact if and only if so
is $T_{n-1}^i$, but there is not reason to believe that, for not
symmetric multilinear forms, the fact that $T_1^i $ is weakly
compact should imply that $T_1^j$ is weakly compact, too. The main
idea behind our proof is to find a trilinear form such that (using
the notation of \cite{ArGa}) $T_1^1$ is weakly compact but $T_1^2$
is not. In \cite{GoGu}, a trilinear form on $\ell_{\infty}$ is
used, which is a slight modification of a trilinear form defined in
\cite{ArChLl}. This form does exactly what we want it to do.}
 \end{remark}

 \begin{remark}
{\rm Although $\hat{T}^1$ is not completely continuous when
considered as a linear operator, its bilinear counterpart
$$T^1: C(K_1) \times
C(K_2) \mapsto (C(K_2))^*$$ given by
$$ T^1(f,g)(h)=T(f,g,h)$$
\begin{em}is\end{em} completely continuous,
according to the usual definition of completely continuous
\begin{em}bilinear\end{em} operator (i.e., if $(f_1^n)\subset C(K_1)$ and
$(f_2^n)\subset C(K_2)$ are weakly Cauchy sequences, then
$T^1(f_1^n, f_2^n)$ is a norm Cauchy sequence), as follows from
\cite{Pe2} (indeed, it follows from \cite{Vi} that
\begin{em}every\end{em}
bilinear continuous operator from $C(K_1) \times C(K_2)$ into
$(C(K_2))^*$ is completely continuous). This proves a conjecture
of \cite{Vi} that states that the fact that a multilinear operator
from the product of $C(K)$ spaces is completely continuous
(considered as a multilinear mapping) does not imply that the same
operator, when considered as a linear mapping from the projective
tensor product of the spaces, has to be completely continuous. }
\end{remark}

It is well known that, for every Banach spaces $E_1, \ldots, E_n$,
$(n>1)$, $E_1 \hat{\otimes} \cdots \hat{\otimes} E_{n-1}$ is
complemented in $E_1 \hat{\otimes} \cdots \hat{\otimes} E_{n}$.
Using this, the next corollary follows easily.

\begin{corolario}
Let $K_1,\ldots K_n$ be infinite compact Hausdorff spaces. Then
\linebreak $C(K_1)\hat{\otimes}\cdots\hat{\otimes}C(K_n) $ has the
DPP if and only if $K_1$,\ldots,$K_n$ are all scattered.
\end{corolario}

The theorem and corollary above remain true for the
\begin{em}symmetric\end{em} projective tensor product.

\begin{teo}
Let $K$ be a compact Hausdorff space. Then $C(K)\hat{\otimes}_s
C(K)$ has the DPP if and only if $K$ is scattered.
\end{teo}
\begin{proof}
If $K$ is scattered, then $(\hat{\otimes}_{n,s} C(K))^*$ is a Schur
space for every $n\in \Bbb N$, and therefore $C(K)\hat{\otimes}_s
C(K)$ has the DPP. Now, if $K$ is not scattered, we can consider
the trilinear form $$ T: C(K) \times C(K) \times C(K) \mapsto \Bbb
K $$ defined by
$$ T(f,g,h)=\sum_{n=1}^{\infty} \frac{1}{2}\left( \xi_n(f) q(g)_n +
\xi_n(g) q(f)_n \right) q(h)_n \quad,$$ that is, the symmetrized
respect to the two first variables of the trilinear form used in
Theorem \ref{teoprincipal}. Now we can apply analogous reasonings
as before to conclude that the linear operator
    $$\hat{T}^1:C(K) \hat{\otimes} C(K) \mapsto
(C(K))^*$$ defined by
    $$ \hat{T}^1(f,g)(h)=T(f,g,h)$$  is weakly compact and not
completely continuous.
\end{proof}

Again, it is well known (see \cite{Bl}) that, for every Banach
space $E$, $\hat{\otimes}_{n-1,s} E$ is complemented in
$\hat{\otimes}_{n,s} E$. So, the following corollary follows.

\begin{corolario}
Let $K$ be a compact Hausdorff space. Then, for every $n>1$,
$\hat{\otimes}_{n,s}C(K) $ has the DPP if and only if $K$ is
scattered.
\end{corolario}

As suggested by J. Guti'rrez, the next theorem can be proved with
the same proof as above.

\begin{teo}
Let $E $ be a Banach space such that

i) Every linear operator from $E$ into $E^*$ is completely
continuous,

ii) $E$ is not Schur, and

iii) $E$ has an isomorphic copy of $\ell_1$,

then $E\hat{\otimes}E$ and $E\hat{\otimes}_s E$ do not have the
DPP.
\end{teo}

There are several spaces which are not $C(K)$ spaces and verify the
conditions of the proposition above, among then,
the space $H^\infty$
\cite{Bo1,Bo2}, the disc algebra \cite{De,Sa}, the space of
analytic uniformly convergent Fourier series on the unit circle
\cite{Sac}.

The authors want to thank Ver¢nica Dimant and Joaqu¡n Guti'rrez
for some clarifying conversations.

\end{document}